\documentclass{amsart}
\usepackage{latexsym}

\theoremstyle{plain}
\newtheorem{Thm}{Theorem}[section]
\newtheorem{Lma}[Thm]{Lemma}

\newtheorem{Prop}[Thm]{Proposition}

\newtheorem{Remark}[Thm]{Remark}
\theoremstyle{definition}
\newtheorem{Def}[Thm]{Definition}

\theoremstyle{remark}

\newcommand{\op}[1]{\operatorname{#1}}

\newcommand{\mlabel}[1]%
  {\mbox{}\marginpar{\raggedleft\hspace{0pt}{\rm\ttfamily#1}}\label{#1}}

\newcommand{\length}{\operatorname{\lambda}}

\newcommand{\fm}{{\mathfrak m}}

\newcommand{\ringR}{\text{$(R,\fm,k)$ }}

\newcounter{hours}\newcounter{minutes}

\newcommand{\excise}[1]{}
\title{\bf On Rings with small Hilbert--Kunz multiplicity}
\author[M.~Blickle]{Manuel Blickle}
\author[F.~Enescu]{Florian Enescu}
\address{FB6 Mathematik, Universit\"at Essen, 45117 Essen, Germany}
\email{manuel.blickle@uni-essen.de}
\address{Department of Mathematics, University of Utah, Salt Lake City,
UT  84112, USA and the Institute of Mathematics of the Romanian
Academy, Bucharest, Romania}
\email{enescu@math.utah.edu}

\date{}

\begin{document}
\begin{abstract} A result of Watanabe and Yoshida says that an unmixed
  local ring of positive characteristic is regular if and only if its Hilbert-Kunz multiplicity is one. We show that, for fixed $p$ and $d$, there exist a number $\epsilon(d,p) > 0$ such that any nonregular unmixed ring $R$ its Hilbert-Kunz multiplicity is at least $1+\epsilon(d,p)$. We also show that local rings with sufficiently small Hilbert-Kunz multiplicity are Cohen-Macaulay and $F$-rational.

\end{abstract}

\maketitle

\section{Introduction}
By a result of Watanabe and Yoshida \cite{WY},
an unmixed local ring $R$ of characteristic $p>0$ is
regular if and only if the Hilbert-Kunz multiplicity,
\[
    e_{HK}(R)=\lim_{e \to \infty} \frac{\length(R/\fm^{[p^e]})}{p^{de}},
\]
is equal to one. A short proof of this was given by Huneke and Yao
in~\cite{HY}. The Hilbert-Kunz multiplicity has proven to be
difficult to compute and even basic questions, such as whether it
is rational, defied a solution despite intense effort. The above
result raises the following, also very basic, question: What is
the smallest possible Hilbert-Kunz multiplicity of a
\emph{non-regular} unmixed ring, say of fixed dimension $d$ and
characteristic $p$? That is, we ask to determine the number
\[
    \epsilon_{HK}(d,p) = \inf\{e_{HK}(R)-1 : \text{$R$ non--regular,
    unmixed, $\dim R = d$, $\op{char} R = p$} \}.
\]
In this paper we start a first investigation into this question by
showing that, at least, $\epsilon_{HK}(d,p)$ is always
\emph{strictly} positive, i.e\ the Hilbert-Kunz multiplicity of a
non-regular ring of fixed dimension and characteristic, cannot be
arbitrarily close to one. The bound we obtain is by no means
optimal as the few cases, where $\epsilon_{HK}(d,p)$ is known, show.

In the course of our investigation we obtain a result of a similar
flavor. We show that a small Hilbert-Kunz multiplicity, precisely
$e_{HK}(R)\leq 1+1/d!$, implies Cohen-Macaulayness and
$F$-rationality (see Proposition \ref{CM}). This sheds some light
into the difficult task of determining the subtle manner in which
the Hilbert-Kunz multiplicity encodes information about the
singularity of $R$.

\section{Lower bounds on the Hilbert-Kunz Multiplicity}

We fix the notation of a $d$--dimensional local ring $\ringR$ with
maximal ideal $\fm$ and residue field $k$ of characteristic $p>0$.
For every positive integer $e$, let $q=p^e$. If $I$ is an ideal of
$R$, then $I^{[q]}=(i^q : i \in I)$.

\begin{Def}
Let $I$ be an $\fm$-primary ideal in $R$. The \emph{Hilbert-Samuel
multiplicity} of $I$ with respect to $R$ is defined by $\lim_{n
\to \infty} \frac{d! \length(R/I^n)}{n^d}$. This limit exists and
is denoted by $e(I,R)$, or, simply, by $e(I)$. The
\emph{Hilbert-Kunz multiplicity} of $I$ with respect to $R$ is
defined by $\lim_{n \to \infty} \frac{\length(R/I^{[q]})}{q^d}$.
This limit exists and is denoted by $e_{HK}(I,R)$, or, simply,
$e_{HK}(I)$. The number $e(\fm,R)$ is called the Hilbert-Samuel multiplicity
of $R$ and usually denoted by $e(R)$.
\end{Def}

For a collection of basic results regarding these notions, we
refer the reader to \cite{Hu}, or \cite{WY}. We list here only a
few well known facts that will be needed later.

One has $\length(R/\fm ^{[q]}) \geq q^d$ and the equality for one
particular $q$ implies that $R$ is regular and hence
$\length(R/\fm ^{[q]})= q^d$ for all $q$. In general, $\max (1,
e(I)/d!) \leq e_{HK}(I) \leq e(I)$, for every $\fm$-primary ideal
of $R$. If $I$ is generated by a system of parameters, then
$e(I)=e_{HK}(I)$. In \cite{Ha} Hanes recently showed that, in fact,
$ e(I)/d! < e_{HK}(I)$, for every $\fm$-primary ideal $I$,
answering affirmatively a question raised by Watanabe and Yoshida
(Question 2.9 in \cite{WY}).

In our paper, we look at the class of unmixed rings. We say that $R$ is unmixed if all the associated primes have the codimension equal to the dimension of $R$. We restrict ourselves to this case, because there are examples of nonregular rings that are not unmixed with $e_{HK}=1$. We will also assume that $d \geq 2$, since in the one dimensional case the Hilbert-Kunz multiplicity is at least $2$ for nonregular rings.

The following lemma, although easy to prove, turns out to be a
very useful tool in many instances (see~\cite{HY}, Lemma 2.1, or
~\cite{WY}, Lemma 4.2).

\begin{Lma}
\label{lemma}
Let $\ringR$ of characteristic $p>0$.

1) $e_{HK}(I) \leq \length(R/I^\ast) \cdot e_{HK}(R)$.

2) For $I$ an $\fm$-primary ideal of $R$, denote $f_I  =
\length(R/I)e_{HK}(R) -e_{HK}(I)$. Then, for every pair of ideals $J
\subset I$ one has $f_I \leq f_J$.
\end{Lma}

The following theorem was conjectured by Watanabe and Yoshida and
is the key tool in our investigation.

\begin{Thm}[Goto-Nakamura, \cite{GN}]
\label{gn} Let $(R,\fm)$ be a homomorphic image of a
Cohen-Macaulay local ring of characteristic $p>0$.

(1)  Assume that $R$ is equidimensional. Then $e (I) \geq
\length(R/I^\ast)$ for every parameter ideal $I$.

(2)  Assume that $R$ is unmixed. If $e (I) = \length(R/I^\ast)$
for some parameter ideal $I$, then $R$ is a Cohen-Macaulay
$F$-rational local ring.
\end{Thm}

\begin{Remark}
{\rm An alternate proof of this theorem has also been given by
Ciuperc\v{a} and Enescu,~\cite{CE} (Theorem 1.1 and Remark 1.10),
where (2) is obtained under some other mild conditions on $R$.}
\end{Remark}

First we present a result which shows that rings with sufficiently
small Hilbert-Kunz multiplicity are Cohen-Macaulay and
$F$-rational.

\begin{Prop}
\label{CM} Let $R$ be an unmixed ring that is a homomorphic image
of a  Cohen-Macaulay  local ring  of characteristic  $p  > 0$. If
$e_{HK}(R) \leq 1 + \op{max}\{1/d!,1/e(R)\}$, then $R$ is
Cohen-Macaulay and $F$-rational.
\end{Prop}

\begin{proof}
As the associated primes of $R$ and $R[x]_{\fm}$ are in bijection,
$R$ is unmixed if and only if $R[x]_{\fm}$ is unmixed. The same holds for the
multiplicities, Cohen--Macaulayness and $F$--rationality. Thus we
can assume that the residue field of $R$ is infinite. Hence, we
can choose a minimal reduction ideal $I$ for $\fm$
(i.e.~$e(I)=e(\fm)=e(R)$) such that $I$ is generated by a system
of parameters.

Let us assume that $R$ is not Cohen-Macaulay and $F$-rational.
We show that this implies that $e_{HK}(R) > 1 + 1/e(R)$. By
Lemma~\ref{lemma} we have that $e(I) \leq \length
(R/I^\ast)e_{HK}(R)$ and substituting $\length(R/I^{\ast}) \leq
e(R) - 1$ (Theorem~\ref{gn}) we get $e(I) \leq (e(I)-1)e_{HK}(R)$.
Thus
\[
    e_{HK}(R) \geq 1+\frac{1}{e(R)-1} > 1+\frac{1}{e(R)}.
\]
To see that $e_{HK}(R) > 1 + 1/d!$, it is now enough to
consider the case $d! \leq e(R)-1$. Then we have
\[
    e_{HK}(R) > \frac{e(R)}{d!} \geq \frac{d!+1}{d!} = 1 + \frac{1}{d!}
\]
finishing the argument.
\end{proof}

\begin{Remark}
\label{dim2}

{\rm Watanabe and Yoshida (\cite{WY}, \cite{WY2}) have shown that, in
dimension $2$, under the assumption of Cohen-Macaulayness, the minimal
Hilbert-Kunz multiplicity is at least 3/2. Our Proposition~\ref{CM}
shows that in fact this assumption can be dropped in the case of
unmixed rings that are a homomorphic image of a Cohen-Macaulay ring.}

\end{Remark}

Now we can address the existence of a lower bound for the
Hilbert-Kunz multiplicity of nonregular rings.
\begin{Thm}
\label{lbound} If $\ringR$ is an unmixed local ring that is the
homomorphic image of a Cohen-Macaulay ring of
characteristic $p>0$ and dimension $d$. If $R$ is not regular,
then $e_{HK} (R) > 1+ \max\{{1}/{(p^dd!)},{1}/{(p^de(R))}\}$.
\end{Thm}

\begin{proof}

As before we can assume that the residue field is infinite and let
$I$ be a minimal reduction for $\fm$. Proposition~\ref{CM} shows
that it is enough to consider the case when $R$ is Cohen-Macaulay
and $F$-rational.

By Lemma~\ref{lemma} we have $f_{\fm ^{[p]}} \leq f_{I^{[p]}}$.
Substituting into this inequality the identities
$e_{HK}(I^{[p]})=e_{HK}(I)p^d$, $e(I^{[p]})=\length(R/I^{[p]})$
and $e(I^{[p]})=e_{HK}(I^{[p]})$ we get
\begin{equation}\label{eq}
    e_{HK}(R)\cdot(\length(R/\fm^{[p]})-p^d) \leq
    e(I^{[p]})\cdot(e_{HK}-1)=p^de(R)\cdot(e_{HK}-1).
\end{equation}
If $R$ is not regular, then $(\length(R/\fm^{[p]})-p^d) \geq 1$ and
we obtain from (\ref{eq}), by solving for $e_{HK}(R)$, that
\[
    e_{HK}(R) \geq 1+\frac{1}{p^de(R)-1}> 1+\frac{1}{p^de(R)}.
\]
Using the inequality $e(R) < d!e_{HK}(R)$ in the right side of
(\ref{eq}) we obtain similarly that
\[
    e_{HK}(R) > 1 + \frac{1}{p^dd!}.
\]
\end{proof}
\begin{Remark}
{\rm Note that this given bound by no means improves the trivial bound
for rings $R$ with $e(R)>d!$ given by
\[
    e_{HK}(R) > e(R)/d! \geq (d!+1)/d! = 1+1/d!.
\]
It is the more interesting case of comparatively small
multiplicity $e(R)$ in which our bound gives new results.}
\end{Remark}

\begin{excise}{
Our results are obtained by manipulating inequalities involving $e_{HK}(R)$,
$e(R)$ and colength of various $\fm$-primary ideals. Improving these
inequalitites will result improving the bounds on the Hilbert-Kunz
multiplicity. In what follows, we present a situation like this. Hanes
has recently proved that the inequality $$e_{HK} \geq e(R) \cdot
\frac{\nu}{(\nu^{1/(d-1)}- 1)^{d-1}}$$ where $\nu$ is the minimal
number of generators for $\fm$. But $e(R) \geq \nu - d +1$ so in fact
$$e_{HK} \geq \frac{\nu-d+1}{d!} \cdot \frac{\nu}{(\nu^{1/(d-1)}-
1)^{d-1}}.$$

Using this inequality, Hanes obtains lower bounds for the Hilbert-Kunz
multiplicity of nonregular rings. He notices that the function $f(\nu)
= \frac{\nu-d+1}{d!} \cdot \frac{\nu}{(\nu^{1/(d-1)}- 1)^{d-1}}$ is
incresing. So $e_{HK} \geq f(d+1)$ (since $\nu \geq d+1$). For $d \leq
4$, $f(d+1)$ is greater than $1$ and lower bounds for the Hilbert-Kunz
multiplicity are obtained).

Let us show how we can use this circle of ideas to improve our
bounds. Firs note that if $R$ is not Cohen-Macaulay, then $e_{HK} \geq
f(d+2)$ ($\nu = d+2$ implies that $R$ is complete intersection, hence
Cohen-Macaulay).

\begin{Prop}
Let $\ringR$ an unmixed local ring of dimension $d$. If $e_{HK}(R)
\leq 1+f(d+2)$, then $R$ is Cohen-Macaulay and $F$-rational.
\end{Prop}

\begin{proof}
Assume that $e_{HK}(R) \leq 1+f(d+2)$. As previously shown, if $R$ is
not Cohen-Macaulay and $F$-rational, then $e(R) \geq 1/f(d+2)
+1$. Hence, $e_{HK} \geq e(R) f(\nu)/(\nu-d+1) \geq (1+f(d+2))/(\nu -d
+1)$.
\end{proof}
}\end{excise}%

\section{Proposed Problems}
A shortcoming of the bound of Theorem \ref{lbound} is
its dependence upon the characteristic. In fact we suspect that
the number
\[
    \epsilon_{HK}(d) = \inf\{\epsilon_{HK}(d,p) : p > 0\}
\]
is strictly bigger than $0$. In the one dimensional case this is
trivial, $\epsilon_{HK}(1) = 1$. Work of Watanabe and
Yoshida~\cite{WY}, together with our Remark~\ref{dim2}, shows that
$\epsilon_{HK}(2) =1/2$. Furthermore, Watanabe and Yoshida give a classification of all two dimensional
Cohen-Macaulay rings with Hilbert-Kunz multiplicity less than 2.
Up to dimension $4$, Hanes's results, \cite{Ha}, also
give lower bounds for $\epsilon_{HK}(d,p)$ that are independent of $p$ (he  can show that $\epsilon_{HK}(3) \geq 1/3$, and that $\epsilon_{HK}(4,p) \geq 0.16$). In private communications, Watanabe mentioned that he and Yoshida can also prove that  $\epsilon_{HK}(3)=1/3$ in dimension $3$.  In dimensions $2$ and $3$, Buchweitz and Chen have obtained the minimal Hilbert-Kunz multiplicity of homogeneous hypersurfaces of fixed degree (see ~\cite{BC}). It should be noted that none of these techniques are known to produce bounds in the general case.

Another question that comes to mind is whether
$\epsilon_{HK}(d,p)$ is attained. If this is the case, what is the
significance of such rings with minimal Hilbert--Kunz
multiplicity? Are they unique up to isomorphism? Again, the work of Watanabe
and Yoshida gives an answer in dimension two only, and all
remaining dimensions remain completely open.

In view of Proposition \ref{CM} we can shift attention away from
non-regular rings to other classes of singularities, and ask for
the minimal Hilbert-Kunz multiplicity of $d$--dimensional ring with
prescribed singularity. The intuition behind this is that bad
singularities should force a high Hilbert-Kunz multiplicity. Our
Proposition \ref{CM} is a result of this type.

{\bf Acknowledgements} 

The second author would like to thank Craig
Huneke, Kei-ichi Watanabe for useful discussions and Douglas Hanes
for providing him with a preprint of his work and interesting correspondence.
The final version of this paper was written while the authors were
visiting the Mathematical Research Sciences Institute at Berkeley, CA,
and they would like to thank this institution for a great work environment and financial support.

\end{document}